\documentclass[11pt]{article}
\usepackage{epsf,amssymb,latexsym,amsmath}
\usepackage{graphicx}
\usepackage{tikz}
\usetikzlibrary{arrows}

\newcommand{\proof}{\noindent{\bf Proof.\ }}
\newcommand{\qed}{\hfill $\square$ \bigskip}
\newcommand{\qqed}{\hfill $\square$ \medskip}

\newtheorem{theorem}{\bf Theorem}[section]

\newtheorem{lemma}[theorem]{\bf Lemma}
\newtheorem{proposition}[theorem]{\bf Proposition}
\newtheorem{conjecture}[theorem]{\bf Conjecture}

\newtheorem{question}[theorem]{\bf Question}
\newtheorem{problem}[theorem]{\bf Problem}
\newtheorem{observation}[theorem]{\bf Observation}

\begin{document}

\title{Hamiltonicity of Cartesian products of graphs\footnote{This work is supported by ARIS project   BI-ME/23-24-022.}}

\author{
    Irena Hrastnik Ladinek\thanks{University of Maribor, FME, Maribor, Slovenia.} \and 
    \v Zana Kovijani\' c Vuki\' cevi\' c\thanks{University of Montenegro, Faculty of Science and Mathematics, Podgorica, Montenegro.} \and
  Tja\v sa Paj Erker\footnotemark[2] \and
   Simon \v Spacapan\footnotemark[2] \footnote{IMFM, Jadranska 19, 1000 Ljubljana.} \footnote{The author is supported by ARIS program P1-0297 and project   
N1-0218.}}
\date{\today}

\maketitle

\begin{abstract}
 A path factor in a graph $G$ is a factor of $G$ in which every component  is a path on at least two vertices. Let $T\Box P_n$ be the Cartesian product of a tree $T$ and a path on $n$ vertices. Kao and Weng \cite{kao} proved that $T\Box P_n$ is hamiltonian if $T$ has a path factor, $n$ is an even integer and $n\geq 4\Delta (T)-2$. They conjectured that for every $\Delta \geq 3$ there exists a graph $G$ of maximum degree $\Delta$ which has a path factor,  such that for every even  $n< 4\Delta-2$ the product  $G\Box P_n$ is not hamiltonian. In this article we prove this conjecture.

\end{abstract}

\noindent
{\bf Key words}:  Cartesian product, hamiltonian graph.

\bigskip\noindent
{\bf AMS subject classification (2020)}: 05C45,  05C76.

% ===================================================================

\bigskip

\section{Introdution}
A {\em  hamiltonian cycle}  in a graph $G$ is a spanning cycle in $G$. A graph is {\em hamiltonian} if it has a hamiltonian cycle. 
The hamiltonicity of graphs is one of the most studied concepts in graph theory, see survey \cite{survey} for a recent  overview of the problem. 
In this article we study the hamiltonicity of Cartesian products of graphs.   
%A graph $G$ is {\em traceable} if it has a path containing  all vertices of $G$. 

%In this article we investigate hamiltonicity of  Cartesian products of two graphs 
%in which one of the factors is a path. 

Let $T$ be a tree and $C_n$ a cycle on $n$ vertices. In \cite{batagelj} the authors proved that the Cartesian product $T\Box C_n$ is hamiltonian if and only if $\Delta(T)\leq n$, where $\Delta(T)$ is the maximum degree of $T$.  
For a graph $G$ let  $\omega(G)$ be the number of components of $G$ and let $t$ be a positive number.  A graph $G$ is {\em $t$-tough} if 
$|S|\geq t \omega(G-S)$ for every separating set $S$ 
of $G$. It is well known that every hamiltonian graph is 1-tough but the converse is not true, see \cite{bauer}. The famous conjecture of Chv\'atal \cite{chvatal} asserts that every  graph which is tough enough is hamiltonian.
A simple corollary of the result obtained in \cite{batagelj}  is that 
$T\Box C_n$ is  hamiltonian if and only if $T\Box C_n$ is $1$-tough.  For products of trees with paths a similar equivalence is not true in general 
 as  there exists  a tree $T$ and a path $P$  such that  $T\Box P$ is $1$-tough but not hamiltonian (see \cite{kao} for details). 
However, in \cite{kao}  the authors  give conditions  under which hamiltonicity of $T\Box P_n$ is equivalent to  $T\Box P_n$ being $1$-tough.  They prove the following theorem. 

\begin{theorem}
Let $T$ be a tree and $n$ an integer. If $T$ has a perfect matching or $n$ is an even integer greater or equal to $4\Delta(T)-2$ then 
 $T\Box P_n$ is hamiltonian if and only if $T\Box P_n$ is 1-tough.  
\end{theorem}
The above theorem is a corollary of the following theorem (also proved in \cite{kao}). 

\begin{theorem} \label{bib}
For any tree  $T$  the Cartesian product $T \Box P_n$ is hamiltonian if one of the following holds:
\begin{itemize}
\item[(a)] $T$ has a perfect matching and $n\geq \Delta(T)$
\item[(b)] $T$ has a path factor, $n$ is an even integer, and $n\geq 4\Delta (T)-2$.
\end{itemize}
\end{theorem}

Note that in (a) and (b)  above, %(a) and (b),
 $T$ is a tree  that has a path factor. 
 We prove that $G\Box H$ is not hamiltonian if $G$ has no path factor and $H$ has a vertex of degree one 
(see  Proposition \ref{pob}). 
It follows that  $T\Box P_n$ is not hamiltonian if  $T$ is a tree with no path factor. Hence, when the hamiltonicity of $T\Box P_n$ is in question,  we can restrict ourselves to trees $T$ that have a path factor. The main result of this article is a proof of the following  conjecture given in \cite{kao}.  

\begin{conjecture} \label{hinezi}
For $k\geq 3$, there is a connected graph $G$  with a path factor such that $\Delta (G)=k$ and $P_{4k-4}\Box G$ is not hamiltonian.  
%Here we prove this conjecture,  which means that the lower bound $n\geq 4\Delta (G)-2$ from in Theorem \ref{bib}  is sharp.
\end{conjecture}
Note that the positive solution to the above conjecture implies that the   bound $n\geq 4\Delta (T)-2$   in Theorem \ref{bib} (b) is sharp. 

In the sequal we list few other known results on the hamiltonicity of Cartesian products of graphs.  
In \cite{cada} authors study graphs that have a 2-factor. If  $F$ is a graph that has a 2-factor then let $g'(F)$ be the minimum length of a  cycle in a 2-factor of $F$  (where the minimum is taken over all 2-factors of $F$). 
If $G$ and $H$ are connected graphs such that both have a  2-factor and $\Delta(G)\leq g'(H)$, 
$\Delta(H)\leq g'(G)$,  then $G\Box H$ is hamiltonian (see  \cite{cada}).  
In \cite{zaks} the author proved that if $n$ is a positive integer and $G^n$ is hamiltonian,
where $G^n$ is the $n$-th Cartesian power of $G$, then $G^k$ is hamiltonian  for every $k\geq n$. 
Hamiltonicity of Cartesian products was also studied in   \cite{dima}, \cite{gauci} and \cite{wallis}.

In the rest of the introduction we give the notation and the terminology and we also prove Proposition \ref{pob}.
For any positive integer $n$ we define $[n]=\{1,\ldots,n\}$. Let $P_n$ be a path on $n$ vertices and let $V(P_n)=[n]$ (we use this notation throughout the paper). A set of pairwise disjoint paths $\mathcal P$ in a graph $G$ is called a {\em path factor} of $G$ if each path  
of $\mathcal P$ has at least two vertices and $\displaystyle\cup_ {P\in \mathcal P} V(P)=V(G)$.   If  $\mathcal P$ is a path factor of a graph $G$, then 
we denote by $\mathcal {E(P)}$ the set of endvertices of paths in  $\mathcal P$. The following observations are straightforward to prove. 

\begin{observation}\label{ena}
If  $\mathcal P$  is a path factor of a graph $G$, then for every set  $\mathcal N\subseteq  \mathcal P$, every vertex of degree 1 in $G-\displaystyle\bigcup_ {P\in \mathcal N} V(P)$ is contained in $\mathcal {E(P)}$.
\end{observation}

\begin{observation}\label{dva}
If  $\mathcal P$  is a path factor of a graph $G$, then for every set  $\mathcal N\subseteq  \mathcal P$, every component of $G-\displaystyle\bigcup_ {P\in \mathcal N} V(P)$ contains a positive even number of vertices in $\mathcal {E(P)}$.
\end{observation}

Let $G=(V(G),E(G))$ and $H=(V(H),E(H))$ be graphs. The {\em Cartesian product} of graphs $G$ and $H$ is the graph, denoted as 
$G\Box H$, with vertex set $V(G\Box H)=V(G)\times V(H)$ where 
vertices $(g_1,h_1)$ and   $(g_2,h_2)$ are adjacent in $G\Box H$ if $g_1=g_2$ and $h_1h_2\in E(H)$ or 
$g_1g_2\in E(G)$ and $h_1=h_2$.

For a graph $G$ we denote by $i(G)$ the number of isolated vertices of $G$.  In \cite{aki} the following criteria for the existence of a path factor is given. 

\begin{proposition}
A graph $G$ has a path factor if and only if for every $S\subseteq V(G)$ we have $i(G-S)\leq 2|S|$. 
\end{proposition}
 Theorem 2 given in \cite{kao} asserts that $G\Box P_n$ is not 
$1$-tough (and hence not hamiltonian) if $G$ is a bipartite graph without a path factor. We give the following strengthening of this result. 

\begin{proposition} \label{pob}
Let $G$ be a graph with no path factor and let $H$ be a graph  with a vertex of degree one. Then 
$G\Box H$ is not hamiltonian. 
\end{proposition}

\proof 
Since $G$ has no path factor  there exists a set $S\subseteq V(G)$ such that $i(G-S)>2|S|$. Let $A$ be the set of isolated vertices in $G-S$ and let 
$u$ be a vertex degree one in $H$. Clearly, no vertex in $A$ is adjacent to a vertex in $G-S$. 

Suppose (reductio ad absurdum) that there is a hamiltonian cycle $C$ in $G\Box H$. Define $A_u=A\times \{u\}$ and $S_u=S\times \{u\}$. 
Since $u$ is a vertex of degree one in $H$, we find that in cycle $C$ each vertex in $A_u$ is adjacent to a vertex in $S_u$. Since 
$|A_u|>2|S|$ we find that a vertex in $S_u$ has a degree at least 3 in $C$, a contradiction which proves the proposition. \qed

\section{Proof of Conjecture \ref{hinezi}}

In this section we prove Conjecture \ref{hinezi} by constructing (for each $k\geq 3$) a tree $T$ of maximum degree $k$ such that 
$T\Box P_n$ is not hamiltonian whenever $n<4k-2$ and $T$ has a path factor.

\begin{lemma} \label{prvasoda}
%Let $P_3=x_1,x_2,x_3$ and $P_n=y_1,\ldots, y_n$ be two paths. %, where $n\geq 3$. 
There exists no path cover $\mathcal P$ of $P_3\Box P_n$,  such that $(2,k)\in \mathcal {E(P)}$  for some odd $k$, and 
$\mathcal {E(P)}\subseteq \{2\}\times ([n]\setminus [k-1])$. 
\end{lemma}

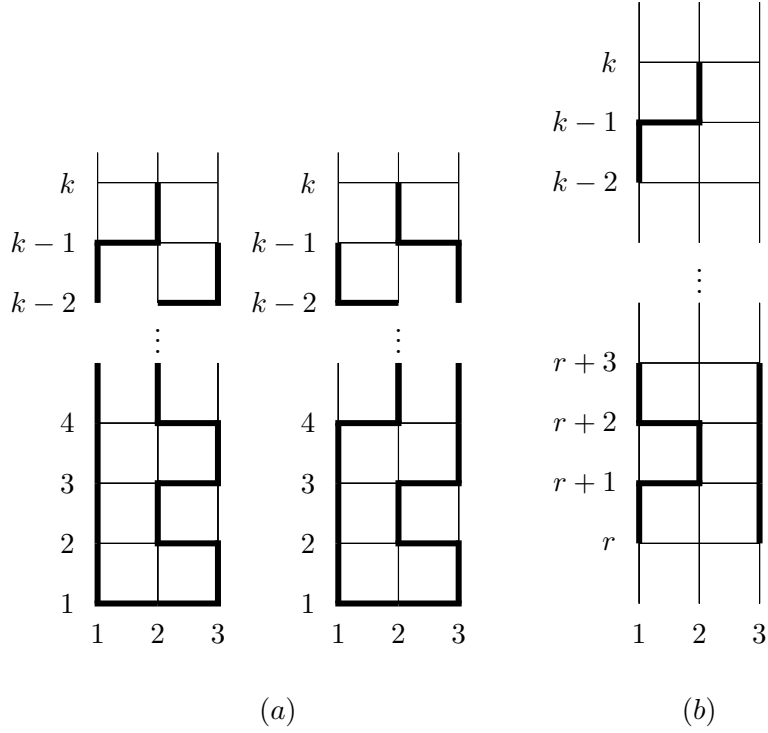
\begin{figure}[htb!]
\begin{center}
\begin{tikzpicture} [scale=0.8]

%prva možnost
\path node at (8,-0.8) {$(a)$};
\path node at (5,0.5) {$1$};
\path node at (6,0.5) {$2$};
\path node at (7,0.5) {$3$};
\path node at (4.5,1) {$1$};
\path node at (4.5,2) {$2$};
\path node at (4.5,3) {$3$};
\path node at (4.5,4) {$4$};
\path node at (4.5,8) {$k$};
\path node at (4.1,7) {$k-1$};
\path node at (4.1,6) {$k-2$};

        \foreach \x in {5,6,7}
	\foreach \y in {1,2,3,4}
	{%navpi�ni
       \draw  (\x,5)--(\x,4);
	\draw (\x,2)--(\x,1);
	\draw (\x,3)--(\x,2);
	\draw (\x,4)--(\x,3);
	%vodoravni
	\draw (5,\y)--(6,\y);
	\draw (6,\y)--(7,\y);
	}
\path node at (6,5.5) {$\vdots$};
	\foreach \x in {5,6,7}
	\foreach \y in {7,8}
	{% navpi�ni
	\draw (\x,7)--(\x,8);
	\draw  (\x,6)--(\x,7);
	\draw  (\x,8)--(\x,8.5);
	% vodoravni
	\draw (5,\y)--(6,\y);
	\draw (6,\y)--(7,\y);
	}
%krepko

\draw [line width=0.8mm] (4.95,1)--(6,1);
\draw [line width=0.8mm] (6,1)--(7.05,1);
\draw [line width=0.8mm] (7,1)--(7,2.05);
\draw [line width=0.8mm] (7,2)--(5.95,2);
\draw [line width=0.8mm] (6,2)--(6,3.05);
\draw [line width=0.8mm] (6,3)--(7.05,3);
\draw [line width=0.8mm] (7,3)--(7,4.05);
\draw [line width=0.8mm] (7,4)--(5.95,4);
\draw [line width=0.8mm] (6,4)--(6,5);

\draw [line width=0.8mm] (5,1)--(5,2);
\draw [line width=0.8mm] (5,2)--(5,3);
\draw [line width=0.8mm] (5,3)--(5,4);
\draw [line width=0.8mm] (5,4)--(5,5);

\draw [line width=0.8mm] (5,6)--(5,7);
\draw [line width=0.8mm] (4.95,7)--(6.05,7);
\draw [line width=0.8mm] (6,7)--(6,8);
\draw [line width=0.8mm] (6,6)--(7.05,6);
\draw [line width=0.8mm] (7,6)--(7,7);

%druga možnost
%\path node at (6,0.5) {$(a)$};
\path node at (9,0.5) {$1$};
\path node at (10,0.5) {$2$};
\path node at (11,0.5) {$3$};
\path node at (8.5,1) {$1$};
\path node at (8.5,2) {$2$};
\path node at (8.5,3) {$3$};
\path node at (8.5,4) {$4$};
\path node at (8.5,8) {$k$};
\path node at (8.1,7) {$k-1$};
\path node at (8.1,6) {$k-2$};

        \foreach \x in {9,10,11}
	\foreach \y in {1,2,3,4}
	{%navpi�ni
       \draw  (\x,5)--(\x,4);
	\draw (\x,2)--(\x,1);
	\draw (\x,3)--(\x,2);
	\draw (\x,4)--(\x,3);
	%vodoravni
	\draw (9,\y)--(10,\y);
	\draw (10,\y)--(11,\y);
	}
\path node at (10,5.5) {$\vdots$};
	\foreach \x in {9,10,11}
	\foreach \y in {7,8}
	{% navpi�ni
	\draw (\x,7)--(\x,8);
	\draw  (\x,6)--(\x,7);
	\draw  (\x,8)--(\x,8.5);
	% vodoravni
	\draw (9,\y)--(10,\y);
	\draw (10,\y)--(11,\y);
	}
%krepko

\draw [line width=0.8mm] (8.95,1)--(10,1);
\draw [line width=0.8mm] (10,1)--(11.05,1);
\draw [line width=0.8mm] (11,1)--(11,2.05);
\draw [line width=0.8mm] (11,2)--(9.95,2);
\draw [line width=0.8mm] (10,2)--(10,3.05);
\draw [line width=0.8mm] (10,3)--(11.05,3);
\draw [line width=0.8mm] (11,3)--(11,4);
\draw [line width=0.8mm] (11,4)--(11,5);
\draw [line width=0.8mm] (10,4)--(10,5);

\draw [line width=0.8mm] (9,1)--(9,2);
\draw [line width=0.8mm] (9,2)--(9,3);
\draw [line width=0.8mm] (9,3)--(9,4);
\draw [line width=0.8mm] (8.95,4)--(10.05,4);
\draw [line width=0.8mm] (10,4)--(10,5);

\draw [line width=0.8mm] (9,7)--(9,6);
\draw [line width=0.8mm] (8.95,6)--(10,6);
\draw [line width=0.8mm] (10,8)--(10,7);
\draw [line width=0.8mm] (9.95,7)--(11.05,7);
\draw [line width=0.8mm] (11,7)--(11,6);
%\end{tikzpicture}
%\caption{??}\label{slika1}

%slika 2
\path node at (15,-0.8) {$(b)$};
\path node at (14,0.5) {$1$};
\path node at (15,0.5) {$2$};
\path node at (16,0.5) {$3$};

\path node at (13.5,2) {$r$};
\path node at (13.1,3) {$r+1$};
\path node at (13.1,4) {$r+2$};
\path node at (13.1,5) {$r+3$};
\path node at (13.5,10) {$k$};
\path node at (13.1,9) {$k-1$};
\path node at (13.1,8) {$k-2$};

        \foreach \x in {14,15,16}
	\foreach \y in {2,3,4,5}
	{%navpi�ni
       \draw  (\x,5)--(\x,4);
	\draw (\x,2)--(\x,1);
	\draw (\x,3)--(\x,2);
	\draw (\x,4)--(\x,3);
            \draw (\x,6)--(\x,5);

	%vodoravni
	\draw (14,\y)--(15,\y);
	\draw (15,\y)--(16,\y);
	}
\path node at (15,6.5) {$\vdots$};
	\foreach \x in {14,15,16}
	\foreach \y in {8,9,10}
	{% navpi�ni
	\draw (\x,7)--(\x,8);
	\draw  (\x,10)--(\x,11);
	\draw  (\x,8)--(\x,9);
          \draw  (\x,9)--(\x,10);

	% vodoravni
	\draw (14,\y)--(15,\y);
	\draw (15,\y)--(16,\y);
	}
\draw [line width=0.8mm] (14,2)--(14,3);
\draw [line width=0.8mm] (13.95,3)--(15.05,3);
\draw [line width=0.8mm] (13.95,4)--(15.05,4);
\draw [line width=0.8mm] (15,3)--(15,3);
\draw [line width=0.8mm] (15,3)--(15,4);
\draw [line width=0.8mm] (15,4)--(15,4);
\draw [line width=0.8mm] (14,4)--(14,5);
\draw [line width=0.8mm] (16,2)--(16,3);
\draw [line width=0.8mm] (16,3)--(16,4);
\draw [line width=0.8mm] (16,4)--(16,5);

\draw [line width=0.8mm] (14,8)--(14,9);
\draw [line width=0.8mm] (13.95,9)--(15.05,9);
\draw [line width=0.8mm] (15,9)--(15,10);
\end{tikzpicture}
\caption{(a) In the proof of Lemma \ref{prvasoda}, the path $P$ contains one of the two subpaths shown on the left. (b) The edges of $E$ in case 2 of Lemma \ref{kukavica}.}
\label{slika1}
\end{center}

\end{figure}

\proof Suppose to the contrary, that there is such a path cover  $\mathcal P$ of $G=P_3\Box P_n$.  Since $(1,1), (3,1)\notin  \mathcal {E(P)}$ are vertices of degree 2 in $G$, we find that there is a path 
 $P\in \mathcal P$ containing subpath  $(1,2),(1,1),(2,1),(3,1),(3,2)$. The vertex $(2,2)\notin  \mathcal {E(P)}$ is also contained in $P$ and  there are only two possibilities: either $(1,2),(2,2),(2,3)$ is a subpath of $P$ or   $(3,2),(2,2),(2,3)$ is a subpath of $P$. By symmetry we may assume the latter, i.e. $P$ contains the subpath  $$(1,2),(1,1),(2,1),(3,1),(3,2),(2,2),(2,3).$$
Now observing vertex $(3,3)\notin  \mathcal {E(P)}$ it has only two  available neighbors (here ‘‘available’’ means possible neighbors of $(3,3)$ in a path 
 $P\in \mathcal P$ containing $(3,3)$), these are $(2,3)$ and $(3,4)$. Hence we found that $(2,3),(3,3),(3,4)$ is a subpath of $P$ and this forces 
that $(1,1),(1,2),(1,3),(1,4)$ is a subpath of $P$. This brings us to the conclusion that 
 $$(1,4),(1,3),(1,2),(1,1),(2,1),(3,1),(3,2),(2,2),(2,3),(3,3),(3,4)$$
is a subpath of $P$. Now we can continue with the same argument (observing the vertex $(2,4)$ forces one of the two symmetric cases: $(1,4),(2,4),(2,5)$ 
or  $(3,4),(2,4),(2,5)$ is a subpath of $P$) which in turn forces the situation shown in Figure \ref{slika1} (a), until eventually we find that either 
$(2,k),(2,k-1),(1,k-1),(1,k-2)$ or $(2,k),(2,k-1),(3,k-1),(3,k-2)$ is a subpath of $P$. In either case 
$(1,k)$ or $(3,k)$ is a vertex of degree 1 in $G-V(P)$, contradicting Observation \ref{ena}. \qed

%\begin{figure}[htb!]
%\begin{center}
%\begin{tikzpicture} [scale=0.8]
%
%%slika 2
%
%\path node at (4.6,1) {$r$};
%\path node at (4.3,2) {$r+1$};
%\path node at (4.3,3) {$r+2$};
%\path node at (4.3,4) {$r+3$};
%\path node at (4.6,9) {$k$};
%\path node at (4.3,8) {$k-1$};
%\path node at (4.3,7) {$k-2$};
%
        %\foreach \x in {5,6,7}
	%\foreach \y in {1,2,3,4}
	%{%navpi�ni
       %\draw  (\x,5)--(\x,4);
	%\draw (\x,2)--(\x,1);
	%\draw (\x,3)--(\x,2);
	%\draw (\x,4)--(\x,3);
            %\draw (\x,0)--(\x,1);
%
	%%vodoravni
	%\draw (5,\y)--(6,\y);
	%\draw (6,\y)--(7,\y);
	%}
%\path node at (6,5.5) {$\vdots$};
	%\foreach \x in {5,6,7}
	%\foreach \y in {7,8,9}
	%{% navpi�ni
	%\draw (\x,7)--(\x,8);
	%\draw  (\x,6)--(\x,7);
	%\draw  (\x,8)--(\x,9);
          %\draw  (\x,9)--(\x,10);
%
	%% vodoravni
	%\draw (5,\y)--(6,\y);
	%\draw (6,\y)--(7,\y);
	%}
%\draw [line width=0.8mm] (5,1)--(5,2);
%\draw [line width=0.8mm] (5,2)--(6,2);
%\draw [line width=0.8mm] (6,2)--(6,3);
%\draw [line width=0.8mm] (6,3)--(5,3);
%\draw [line width=0.8mm] (5,3)--(5,4);
%\draw [line width=0.8mm] (7,1)--(7,2);
%\draw [line width=0.8mm] (7,2)--(7,3);
%\draw [line width=0.8mm] (7,3)--(7,4);
%
%\draw [line width=0.8mm] (5,7)--(5,8);
%\draw [line width=0.8mm] (5,8)--(6,8);
%\draw [line width=0.8mm] (6,8)--(6,9);
%\end{tikzpicture}
%\caption{A step in the proof of Lemma \ref{kukavica}}\label{slika2}
%\end{center}
%\end{figure}
%
%
%
\begin{lemma}\label{kukavica}
%Let $P_3=y_1,y_2,y_3$ and $P_n=x_1,\ldots, x_n$ be two paths. %, where $n\geq 3$. 
There exists no path cover $\mathcal P$ of $P_3\Box P_n$ with the following properties.
\begin{itemize}
\item[(i)] $\mathcal {E(P)}\subseteq \{2\}\times [n]$;
\item[(ii)] $(2,k)\in \mathcal {E(P)}$ for some odd integer $k$ and 
$(2,k-1)\notin \mathcal {E(P)}$, and
\item [(iii)] if $i<k$ and $i$ is odd, then $(2,i)\in \mathcal {E(P)}$ if and only if $(2,i-1)\in \mathcal {E(P)}$.
\end{itemize}
\end{lemma}

\proof Suppose to the contrary that a path cover  $\mathcal P$ with properties 
(i), (ii) and (iii) does exist. 
%Let us first exclude the possibility that a path $P\in \mathcal P$ has endvertices  $(y_2,x_{a})$ and  $(y_2,x_{b})$ such that $a<k$ and $b\geq k$. 
Let $r$ be the maximum integer $i<k$ such that $(2,i)\in \mathcal {E(P)}$ (by Lemma \ref{prvasoda}  $r$ is well defined, and by (ii) and  (iii) $r$ is odd).
Define $E=\bigcup_{P\in \mathcal P}E(P)$. Let $e_1=(1,r)(1,r+1)$ and $e_3=(3,r)(3,r+1)$. Observe that $e_1\notin E$ implies 
 $(1,r)(2,r)\in E$ because $(1,r)\notin  \mathcal {E(P)}$ is a vertex of degree 3 in $G=P_3\Box P_n$. Similarly 
$e_3\notin E$ implies 
 $(3,r)(2,r)\in E$. Since $(2,r)\in  \mathcal {E(P)}$ we find that 
$e_1\in E$ or  $e_3\in E$. 

{\em Case 1:} Suppose that $e_1\notin E$ (the case  $e_3\notin E$ is symmetric). Then 
$f=(1,r)(2,r)\in E$ and $g=(1,r)(1,r-1)\in E$. Assume that $P_1,P_2\in \mathcal P$ are paths such that $f,g\in P_1$ and $e_3\in P_2$ 
(here it is possible that $P_1=P_2$, but the arguments given below apply in either case). % in which case the proof is the same). 
Then $V(P_1)\cup V(P_2)$ contains vertices $(1,r),(2,r)$ and $(3,r)$ and therefore $V(P_1)\cup V(P_2)$ is a separating set in 
$G$. Moreover, by (iii) components of $G-(V(P_1)\cup V(P_2))$ which are contained in $[3]\times [r-1]$ have an odd number of vertices in $\mathcal {E(P)}$. So there exists a component of $G-(V(P_1)\cup V(P_2))$ which has an odd number of vertices in $\mathcal {E(P)}$,
contradicting Observation \ref{dva}. 

{\em Case 2:} Suppose that $e_1,e_3\in E$. If also $e_2=(2,r)(2,r+1)\in E$ then the arguments are similar as in Case 1. Assume that 
$e_1\in P_1,e_2\in P_2$ and $e_3\in P_3$ (here it is possible that $P_1,P_2$ and $P_3$ are not pairwise distinct paths, but the arguments apply in either case). 
Then $S=V(P_1)\cup V(P_2)\cup V(P_3)$ is a separating set  in $G$ such that $G-S$ has a connected component which contains an odd number of vertices in 
 $\mathcal {E(P)}$ contradicting Observation \ref{dva}. 

Thus we may assume that  $e_1,e_3\in E$ and $e_2\notin E$. %Let $e_1\in P_1$ and $e_3\in P_2$ (possibly $P_1=P_2$).
%If $P_1$ contains a vertex in $X=\{3\}\times ([k]\setminus [r])$ then $P_1=P_2$ for otherwise (since $e_2\notin E$)
 % an endvertex of $P_1$ (or $P_2$) is not in $\mathcal {E(P)}$. Let us assume that $P_1$ contains a vertex in $X$ and thus $P_1=P_2$.
If $(1,r+1)(2,r+1)\in E$ and $(2,r+1)(3,r+1)\in E$, then $Y=[3]\times ([n]\setminus [r+1])$ is covered by paths in a set $\mathcal N\subset \mathcal P$. 
Since $k$ and $r$ are odd, we find that $\mathcal N$ is a path cover of $Y$, contradicting Lemma \ref{prvasoda}.

Otherwise (since $e_2\notin E$) $(2,r+1)(2,r+2)\in E$ and either $(1,r+1)(2,r+1)\in E$ or $(2,r+1)(3,r+1)\in E$. By symmetry we may assume 
 $(1,r+1)(2,r+1)\in E$. Observing vertex  $(1,r+2)$, which has only two available neighbors, we find that 
 $(1,r+2)(2,r+2)\in E$ and  $(1,r+2)(1,r+3)\in E$. It follows also that  $(3,r+1)(3,r+2)\in E$ and  $(3,r+2)(3,r+3)\in E$, see Figure \ref{slika1} (b). 
% illustrates this situation.  
A similar argument is repeated until we find that  $(1,k-2),(1,k-1),(2,k-1),(2,k)$ or $(3,k-2),(3,k-1),(2,k-1),(2,k)$ is a subpath of a path $P\in \mathcal P$. It follows that $(1,k)$ or $(3,k)$ is a vertex of degree 1 
in $G-V(P)$, contradicting Observation \ref{ena}.  \qed

\begin{figure}[htb!]
\begin{center}
\begin{tikzpicture} [scale=0.8]

	{    	\draw [line width=0.8mm] (5,5)--(9,5);
			\draw (9,5)--(9,6.5);
			\draw (9,5)--(10.3,4.25);
			\draw (9,5)--(9,3.5);
			
			\draw [line width=0.8mm](9,6.5)--(8.37,7.86);
			\draw [line width=0.8mm](9,6.5)--(9.63,7.86);
			\draw [line width=0.8mm] (10.3,4.25)--(11.16,3.02);
			\draw [line width=0.8mm](10.3,4.25)--(11.79,4.12);
			\draw [line width=0.8mm](9,3.5)--(9.63,2.14);
			\draw [line width=0.8mm] (9,3.5)--(8.37,2.14);
			
%\draw [line width=0.8mm] (2,1)--(3,1);
\filldraw (5,5) circle (3pt);
\path node at (5.3,4.6) {$a$};
\filldraw (7,5) circle (3pt);
\path node at (7,4.6) {$b$};
\filldraw (9,5) circle (3pt);
\path node at (8.7,4.6) {$c$};

\filldraw (9,6.5) circle (3pt);
\path node at (8.3,6.5) {$c_{\Delta -1}$};
\filldraw (10.3,4.25) circle (3pt);
\path node at (10.3,4.6) {$c_2$};
\filldraw (9,3.5) circle (3pt);
\path node at (8.3,3.5) {$c_{1}$};

\filldraw (5,6.5) circle (3pt);
\path node at (5.8,6.5) {$a_{\Delta -1}$};
\filldraw (3.7,4.25) circle (3pt);
\path node at (3.7,4.6) {$a_2$};
\filldraw (5,3.5) circle (3pt);
\path node at (5.8,3.5) {$a_{1}$};

\draw (5,5)--(5,6.5);
\draw (5,5)--(3.7,4.25);
\draw (5,5)--(5,3.5);

\filldraw (8.37,7.86) circle (3pt);
\filldraw (9.63,7.86) circle (3pt);
\path node at (8.37,8.3) {$z_{\Delta -1}$};
\path node at (9.63,8.3) {$y_{\Delta -1}$};

\filldraw (11.16,3.02) circle (3pt);
\filldraw (11.79,4.12) circle (3pt);
\path node at (11.6,2.7) {$y_2$};
\path node at (12.3,4.15) {$z_2$};

\filldraw (9.63,2.14) circle (3pt);
\filldraw (8.37,2.14) circle (3pt);
\path node at (9.7,1.8) {$z_{1}$};
\path node at (8.2,1.8) {$y_{1}$};

\filldraw (5.63,7.86) circle (3pt);
\filldraw (4.37,7.86) circle (3pt);
\path node at (5.63,8.3) {$v_{\Delta -1}$};
\path node at (4.37,8.3) {$u_{\Delta -1}$};

\filldraw (2.84,3.02) circle (3pt);
\filldraw (2.21,4.12) circle (3pt);
\path node at (2.5,2.7) {$u_2$};
\path node at (1.7,4.15) {$v_2$};

\filldraw (4.37,2.14) circle (3pt);
\filldraw (5.63,2.14) circle (3pt);
\path node at (4.3,1.8) {$v_{1}$};
\path node at (5.67,1.8) {$u_{1}$};

\draw [line width=0.8mm](5,6.5)--(5.63,7.86);
\draw [line width=0.8mm](5,6.5)--(4.37,7.86);
\draw [line width=0.8mm](3.7,4.25)--(2.84,3.02);
\draw [line width=0.8mm](3.7,4.25)--(2.21,4.12);
\draw [line width=0.8mm](5,3.5)--(4.37,2.14);
\draw [line width=0.8mm](5,3.5)--(5.63,2.14);

\filldraw (11,6.1) circle (1.3pt);
\filldraw (11.15,5.75) circle (1.3pt);
\filldraw (10.74,6.38) circle (1.3pt);
\filldraw (3,6.1) circle (1.3pt);
\filldraw (2.85,5.75) circle (1.3pt);
\filldraw (3.24,6.38) circle (1.3pt);

\path node at (9,9.1) {$C_{\Delta -1}$};
\path node at (12.3,3.3) {$C_2$};
\path node at (9,1) {$C_{1}$};
\path node at (5,9.1) {$A_{\Delta -1}$};
\path node at (1.6,3.3) {$A_2$};
\path node at (5,1) {$A_{1}$};
\path node at (7,5.5) {$B$};

	}	
\end{tikzpicture}
\caption{The tree $T_\Delta$. A path factor of $T_\Delta$ is denoted by bold lines.}
\label{slika2}
\end{center}
\end{figure}
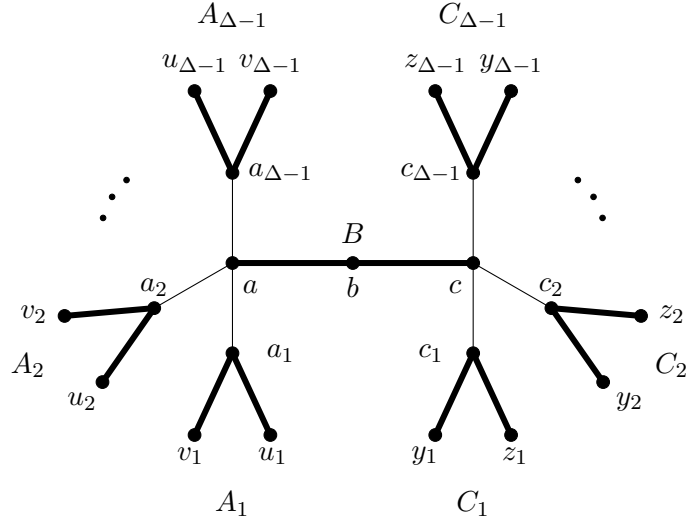

\begin{figure}[htb!]
\begin{center}
\begin{tikzpicture} [scale=0.8]

%Mi Ni...
% ii

        \foreach \x in {5,6,7}
	\foreach \y in {2,3,4,5}
	{%navpi�ni
       \draw  (\x,5)--(\x,4);
	\draw (\x,2)--(\x,1);
	\draw (\x,3)--(\x,2);
	\draw (\x,4)--(\x,3);
          \draw (\x,5)--(\x,6);
	%vodoravni
	\draw (5,\y)--(6,\y);
	\draw (6,\y)--(7,\y);
\draw [line width=0.8mm] (5,2)--(5,3.05);
\draw [line width=0.8mm] (5,3)--(6,3);
\draw [line width=0.8mm] (6,2.95)--(6,4.05);
\draw [line width=0.8mm] (6,4)--(4.95,4);
\draw [line width=0.8mm] (5,4)--(5,5);
\draw [line width=0.8mm] (7,2)--(7,3.05);
\draw [line width=0.8mm] (7,3)--(8,3);
\draw [line width=0.8mm] (7,4)--(8,4);
\draw [line width=0.8mm] (7,3.95)--(7,5);

\path node at (5.5,3.4) {$M_i$};
\path node at (7.5,4.4) {$R_{ij}$};
\path node at (7.5,2.4) {$Q_{ij}$};
	}

%  i
\foreach \x in {1,2,3}
	\foreach \y in {2,3,4,5}
	{%navpi�ni
       \draw  (\x,5)--(\x,4);
	\draw (\x,2)--(\x,1);
	\draw (\x,3)--(\x,2);
	\draw (\x,4)--(\x,3);
          \draw (\x,5)--(\x,6);
	%vodoravni
	\draw (1,\y)--(2,\y);
	\draw (2,\y)--(3,\y);
\draw [line width=0.8mm] (1,2)--(1,3);
\draw [line width=0.8mm] (0.95,3)--(2.05,3);
\draw [line width=0.8mm] (2,3)--(2,4);
\draw [line width=0.8mm] (2.05,4)--(0.95,4);
\draw [line width=0.8mm] (1,4)--(1,5);
\draw [line width=0.8mm] (3,2)--(3,3);
\draw [line width=0.8mm] (3,3)--(3,4);
\draw [line width=0.8mm] (3,4)--(3,5);

\path node at (1.5,3.4) {$M_i$};
\path node at (3.5,3.4) {$N_{i}$};
	}

%  iii
\foreach \x in {9,10,11}
	\foreach \y in {2,3,4,5}
	{%navpi�ni
       \draw  (\x,5)--(\x,4);
	\draw (\x,2)--(\x,1);
	\draw (\x,3)--(\x,2);
	\draw (\x,4)--(\x,3);
          \draw (\x,5)--(\x,6);
	%vodoravni
	\draw (9,\y)--(10,\y);
	\draw (10,\y)--(11,\y);
\draw [line width=0.8mm] (9,2)--(9,3);
\draw [line width=0.8mm] (9,3)--(9,4);
\draw [line width=0.8mm] (9,4)--(9,5);
\draw [line width=0.8mm] (11,2)--(11,3);
\draw [line width=0.8mm] (11.05,3)--(9.95,3);
\draw [line width=0.8mm] (10,3)--(10,4);
\draw [line width=0.8mm] (9.95,4)--(11.05,4);
\draw [line width=0.8mm] (11,4)--(11,5);

\path node at (8.5,3.4) {$T_i$};
\path node at (10.5,3.4) {$S_{i}$};
	}

%  iv
\foreach \x in {13,14,15}
	\foreach \y in {2,3,4,5}
	{%navpi�ni
       \draw  (\x,5)--(\x,4);
	\draw (\x,2)--(\x,1);
	\draw (\x,3)--(\x,2);
	\draw (\x,4)--(\x,3);
          \draw (\x,5)--(\x,6);
	%vodoravni
	\draw (13,\y)--(14,\y);
	\draw (14,\y)--(15,\y);
\draw [line width=0.8mm] (12,3)--(13,3);
\draw [line width=0.8mm] (13,3.05)--(13,2);
\draw [line width=0.8mm] (12,4)--(13,4);
\draw [line width=0.8mm] (13,3.95)--(13,5);
\draw [line width=0.8mm] (15,2)--(15,3);
\draw [line width=0.8mm] (15.03,3)--(13.95,3);
\draw [line width=0.8mm] (14,3)--(14,4);
\draw [line width=0.8mm] (13.95,4)--(15.05,4);
\draw [line width=0.8mm] (15,4)--(15,5);

\path node at (14.5,3.4) {$S_i$};
\path node at (12.5,4.4) {$V_{ij}$};
\path node at (12.5,2.4) {$U_{ij}$};
	}

%  v
\foreach \x in {17,18,19}
	\foreach \y in {2,3,4,5}
	{%navpi�ni
       \draw  (\x,5)--(\x,4);
	\draw (\x,2)--(\x,1);
	\draw (\x,3)--(\x,2);
	\draw (\x,4)--(\x,3);
          \draw (\x,5)--(\x,6);
	%vodoravni
	\draw (17,\y)--(18,\y);
	\draw (18,\y)--(19,\y);
\draw [line width=0.8mm] (17,2)--(17,3.05);
\draw [line width=0.8mm] (17,3)--(18,3);
\draw [line width=0.8mm] (18,3)--(19,3);
\draw [line width=0.8mm] (19,3.05)--(19,2);
\draw [line width=0.8mm] (17,5)--(17,3.95);
\draw [line width=0.8mm] (17,4)--(18,4);
\draw [line width=0.8mm] (18,4)--(19,4);
\draw [line width=0.8mm] (19,3.95)--(19,5);

\path node at (18,2.4) {$X_i$};
\path node at (18,4.4) {$Z_i$};
	}
\end{tikzpicture}
\caption{Ten different paths in $T_n\Box P_{m}$}
\label{slika3}
\end{center}
\end{figure}
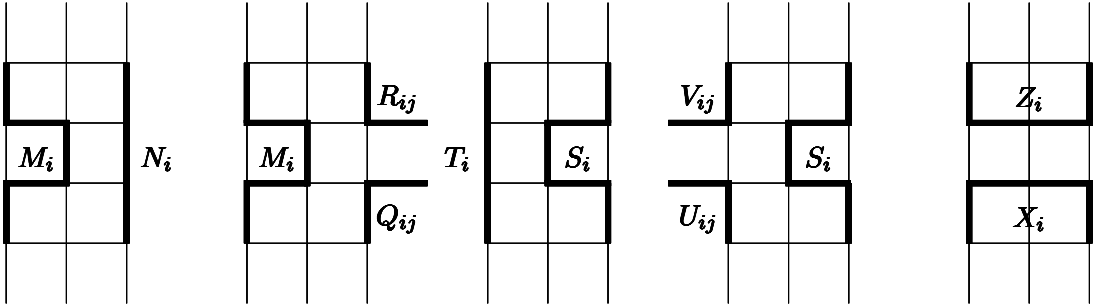

Let $\Delta\geq 3$ be an integer and let $T_\Delta$ be the tree shown in Figure \ref{slika2}, with one vertex of degree 2, two vertices of degree $\Delta$,  $2\Delta-2$ vertices of degree 3, and  
$4\Delta-4$ vertices of degree 1. Let $a$ and $c$ be vertices of degree $\Delta$ and $b$ the vertex adjacent to $a$ and $c$. Every vertex of degree 3   in $T_\Delta$ is adjacent to two vertices of degree 1. Vertices of degree 3 are denoted by 
$a_1,\ldots,a_{\Delta-1}$ and $c_1,\ldots,c_{\Delta-1}$. Vertices of degree 1 adjacent to $a_i$ are $u_i, v_i$ and vertices of degree 1 adjacent to 
$c_i$ are $y_i, z_i$ for $i\in [\Delta-1]$. We define paths in $T_\Delta$: $A_i=u_i,a_i,v_i$ and  $C_i=y_i,c_i,z_i$ for $i\in [\Delta-1]$ and $B=a,b,c$. 
For $i\in [m-3]$ and $j\in [n-1]$ we define the following paths in $T_n\Box P_{m}$ shown in Figure \ref{slika3}:
%
%\begin{itemize}
%\item[(a)] $M_i=(a,x_i),(a,{i+1}),(b,{i+1}),(b,{i+2}),(a,{i+2}),(a,{i+3})$
%\item [(b)] $N_i=(c,i),(c,{i+1}),(c,{i+2})$
%\item[(c)]  $Q_{i,j}=(c,i),(c,{i+1}),(c_j,{i+1})$
%\item[(d)]  $R_{i,j}=(c_j,{i+2}),(c,{i+2}),(c,{i+3})$
%\item[(e)] $S_i=(c,i),(c,{i+1}),(b,{i+1}),(b,{i+2}),(c,{i+2}),(c,{i+3})$
%\item [(f)] $T_i=(a,i),(a,{i+1}),(a,{i+2})$
%\item[(g)]  $U_{i,j}=(a,i),(a,{i+1}),(a_j, {i+1})$
%\item[(h)]  $V_{i,j}=(a_j,{i+2}),(a,{i+2}),(a,{i+3})$
%\end{itemize}

 \begin{quote}
%\smallskip
$M_i=(a,i),(a,{i+1}),(b,{i+1}),(b,{i+2}),(a,{i+2}),(a,{i+3})$\\
\smallskip
 $N_i=(c,i),(c,{i+1}),(c,{i+2}),(c,{i+3})$\\
\smallskip
 $Q_{i,j}=(c,i),(c,{i+1}),(c_j,{i+1})$\\
\smallskip
   $R_{i,j}=(c_j,{i+2}),(c,{i+2}),(c,{i+3})$\\
\smallskip
 $S_i=(c,i),(c,{i+1}),(b,{i+1}),(b,{i+2}),(c,{i+2}),(c,{i+3})$ \\
\smallskip
 $T_i=(a,i),(a,{i+1}),(a,{i+2}),(a,{i+3})$\\
\smallskip
  $U_{i,j}=(a,i),(a,{i+1}),(a_j, {i+1})$\\
\smallskip
  $V_{i,j}=(a_j,{i+2}),(a,{i+2}),(a,{i+3})$\\
\smallskip
 $X_i=(a,i),(a,{i+1}),(b,{i+1}),(c,{i+1}),(c,i)$\\
\smallskip
 $Z_i=(a,i+3),(a,{i+2}),(b,{i+2}),(c,{i+2}),(c,i+3)$ 
\end{quote}

\begin{theorem}
For every $\Delta\geq 3$ there exist a graph $G$ of maximum degree $\Delta$ such that $G$ has a path factor and $G\Box P_{m}$ is not hamiltonian for 
every $m\leq 4\Delta-3$. 
\end{theorem}

\proof
We claim that for every $\Delta\geq 3$ the product $T_\Delta\Box P_{m}$ is not hamiltonian whenever $m\leq 4\Delta-3$. Suppose (reductio ad absurdum) that for some $n\geq 3$ and $m\leq 4n-3$ the graph $H=T_n\Box P_{m}$ is hamiltonian and let  $C$ be a hamiltonian cycle in  $H$.  \\

{\em Claim 1:} For every $i\in [n-1]$, $C$ contains   paths $(u_i,2),(u_i,1),(a_i,1),(v_i,1),(v_i,2)$ and  $(y_i,2),(y_i,1),(c_i,1),(z_i,1),(z_i,2)$.\\
\indent {\em Proof:} The claim follows from the fact that  $(u_i,1),(v_i,1),(y_i,1),(z_i,1)$ are vertices of degree 2 in $H$. \qqed

{\em Claim 2:} $C$ contains the  path $(a,2),(a,1),(b,1),(c,1),(c,2)$.\\
\indent {\em Proof:} By Claim 1, $C$ contains no edge $(a_i,1)(a,1)$ or $(c_i,1)(c,1)$ for $i\in [n-1]$. 
In the remaining graph (once the above mentioned edges are removed from $H$) vertices $(a,1)$ and $(c,1)$ have degree 2. \qqed

{\em Claim 3:} For every odd $i\in [m-3]$, $C$ contains one of the following paths 
\begin{itemize}
\item[(i)] $M_i$ and $N_i$, or 
\item[(ii)] $M_i, Q_{i,j}$ and $R_{i,j}$ for some $j\in [n-1]$, or
\item[(iii)] $S_i$ and $T_i$, or 
\item[(iv)] $S_i, U_{i,j}$ and $V_{i,j}$ for some $j\in [n-1]$, or 
\item[(v)]  $Z_i$ and $X_i$. 
\end{itemize}

Before we prove Claim 3 let us first finish the proof of the theorem (using Claim 3). First observe that there are $2n-2$ components of  $H-V(B\Box P_{m})$.  
By Claim 3 and Claim 2 there are at most $m-2$ edges of $C$ with one endvertex in 
$B\Box P_{m}$ and the other in $H-V(B\Box P_{m})$.  It follows that $C$ intersects at most $(m-2)/2$ components of  $H-V(B\Box P_{m})$ and 
therefore $C$ is not a spanning cycle of $H$ (recall that $m\leq 4n-3$). This contradiction proves the theorem.

\indent {\em Proof of Claim 3:}  We prove Claim 3 by induction on $i$. Let us start with $i=1$. 
By Claim 2, $(a,2),(a,1),(b,1),(c,1),(c,2)$ is a path in $C$. Observing vertex $(b,2)$ we find that either 
 $(a,2),(b,2),(b,3)$  or  $(c,2),(b,2),(b,3)$  is a path in $C$. We assume the former (the proof for the other case is analogous). 
By Lemma \ref{prvasoda}, $C$ does not use any edge $(a,3)(a_j,3)$ for $j\in [n-1]$, and therefore 
 $(b,3),(a,3),(a,4)$ is a path in $C$. 

If $C$ uses an  edge $(c,3)(c_j,3)$ for $j\in [n-1]$ then, by Lemma \ref{prvasoda}, it also uses the edge 
$(c,2)(c_j,2)$ and we have case (ii) of Claim 3 (for $i=1$). Otherwise  $(c,2),(c,3),(c,4)$ is a path in $C$ and we have case (i) of Claim 3  (for $i=1$).

Suppose that the claim is true for all odd $i\leq  k-4$ (where $k-4<m-3$ and $k-4$ is odd). We shall prove it for $i=k-2$. 
First note that for every $j\in [n-1]$ every vertex of  $A_j\Box P_{m}$ has at most one neighbor not in $A_j\Box P_{m}$. 
It follows that components of $C\cap (A_j\Box P_{m})$ 
form a path cover of $A_j\Box P_{m}$. Similarly components of  $C\cap (C_j\Box P_{m})$ form a path cover of $C_j\Box P_{m}$. 
Let  $\mathcal P_j$ be the  path cover obtained from  
$C\cap (A_j\Box P_{m})$ and  $\mathcal R_j$ be the  path cover obtained from  
$C\cap (C_j\Box P_{m})$ for $j\in [n-1]$. Note that  $\mathcal P_j$ and $\mathcal R_j$ have property (i) of Lemma \ref{kukavica}.
Since the claim is true for all odd $i\leq k-4$ we find by (ii) and (iv) that for every $j\in [n-1]$ path covers  $\mathcal P_j$ and $\mathcal R_j$  have %possesses 
the property (iii) of Lemma \ref{kukavica}.

Since the claim is true for $i=k-4$  % either   
 $(a,{k-1}), (a,{k-2})(b,{k-2}),(b,{k-3})$ or  $(c,{k-1}), (c,{k-2})(b,{k-2}),(b,{k-3})$ or  $(a,{k-1}), (a,{k-2})(b,{k-2}),(c,{k-2}),(c,k-1)$  is a path in $C$; moreover  
$(a,{k-1})(a,{k-2})\in E(C)$ and $(c,{k-1})(c,{k-2})\in E(C)$ (see Figure \ref{slika3}). 
Therefore  
(since $(b,{k-1})$ is a vertex of degree 4 in $H$) one of the following   occurs:
\begin{itemize}
\item[I.]  $(a,{k-1})(b,{k-1})\in E(C)$ and  $(b,{k-1})(c,{k-1})\in E(C)$, or
\item[II.]  $(a,{k-1})(b,{k-1})\in E(C)$ and  $(b,{k-1})(b,{k})\in E(C)$, or
\item[III.]  $(c,{k-1})(b,{k-1})\in E(C)$ and  $(b,{k-1})(b,{k})\in E(C)$.
\end{itemize}

{\em Case I.} In this case $(a,k-2),(a,k-1),(b,k-1),(c,k-1),(c,k-2)$  is a path in $C$, and hence $C$ does not use any edge $(a,k-1)(a_j,k-1)$ or 
$(c,k-1)(c_j,k-1)$ for $j\in [n-1]$. 

If $C$ uses edge $(a,k)(a_j,k)$ or $(c,k)(c_j,k)$ for some $j\in [n-1]$, then path covers $\mathcal P_j$ or $\mathcal R_j$  have %possesses 
property (ii) of Lemma \ref{kukavica}, respectively. Recall also that $\mathcal P_j$ and $\mathcal R_j$ have properties (i) and (iii) of Lemma  
 \ref{kukavica}, and therefore, by the same lemma, $\mathcal P_j$ or $\mathcal R_j$ are not path covers of  
$A_j\Box P_{m}$ or $C_j\Box P_{m}$, a contradiction. It follows that $C$ does not use any edge 
$(a,k)(a_j,k)$ or $(c,k)(c_j,k)$ for $j\in [n-1]$, and therefore 
$(a,k+1),(a,k),(b,k),(c,k),(c,k+1)$ is a path in $C$. Hence, we have case (v) of Claim 3  (for $i=k-2$). 

{\em Case II. and Case III.} Both cases are symmetric, so we prove only Case II.
Similarly as in Case I we find that  $C$ does not use any edge 
$(a,k)(a_j,k)$  for $j\in [n-1]$, and therefore 
$(a,k+1),(a,k),(b,k),(b,k-1),(a,k-1),(a,k-2)$ is a path in $C$.

If $C$ uses  edge $(c,k)(c_j,k)$ for some $j\in [n-1]$ then, by Lemma \ref{kukavica}, it has to use also the edge $(c,k-1)(c_j,k-1)$ and we have 
case (ii) of Claim 3 (for $i=k-2$). Otherwise  $(c,k-1),(c,k),(c,k+1)$ is a path in $C$ and we have the case (i) of Claim 3  (for $i=k-2$). 

This proves the claim and hence also the theorem. 
\qed

\section{Concluding remarks}

As mentioned in the introduction, there exists a product $T\Box P_n$, exhibited in \cite{kao}, which is 1-tough but not hamiltonian.
Moreover, in \cite{kral} the authors constructed a sequence of prisms (the products of $G$ with  $K_2$), each of which is not hamiltonian,  and such that the toughness of the prisms approaches $9/4$.
However, we are not aware of any product of 2-connected graphs $G$ and $H$, such that $G\Box H$ is 1-tough but not hamiltonian. Hence we pose the following question.

\begin{question}
If $G$ and $H$ are 2-connected graphs,  is $G\Box H$ hamiltonian if and only if $G\Box H$ is 1-tough ?
\end{question}

If $B$ is a bipartite graph with parts of the bipartition having  equal size  then we say that $B$ is {\em balanced}, otherwise it is {\em unbalanced}.  
If $B$ is unbalanced then it is not 1-tough and hence not hamiltonian. 
However, except the few results given in \cite{kao}, nothing is known regarding toughness of Cartesian products of graphs in general. 
If $B_1$ and $B_2$ are bipartite graphs and none of them is balanced, then it is easy to see that $B_1\Box B_2$ is unbalanced. 
Thus it is required that at least one of $B_1$ and $B_2$ is balanced if we wish $B_1\Box B_2$ to be hamiltonian. However, a full characterization of 
products of bipartite graphs that are 1-tough is still an open question. In particular, we pose the following problem.   

\begin{problem}
Characterize products of trees $T_1\Box T_2$ that are 1-tough.
\end{problem}

%The same problem can be asked for products of bipartite graphs. 
A solution of the above problem would give us a sufficient condition for hamiltonicity of products of trees.

\end{document}